\documentclass{amsart}

\usepackage{fullpage}
\usepackage{epic}

\title{An elliptic hypergeometric integral with $W(F_4)$ symmetry.}
\author{Fokko J. van de Bult \\ California Institute of Technology}
\date{September 18, 2009}

\newtheorem{thm}{Theorem}[section]
\newtheorem{definition}[thm]{Definition}
\newtheorem{cor}[thm]{Corollary}
\newtheorem{prop}[thm]{Proposition}

\newcommand{\rphis}[2]{{}_{#1\vphantom{#2}}\phi_{#2\vphantom{#1}}}

\newcommand{\rWs}[2]{{}_{#1\vphantom{#2}}W_{#2\vphantom{#1}}}

\newcommand{\rphisx}[4]{\rphis{#1}{#2}\left( \begin{array}{c} #3 \end{array};q,#4\right)}

\begin{document}

\begin{abstract}
In this article we give a new transformation between elliptic hypergeometric beta integrals, which
gives rise to a Weyl group symmetry of type $F_4$. The transformation is a generalization of a series
transformation discovered by Langer, Schlosser, and Warnaar \cite{LSW}. Moreover we consider various limits of this
transformation to basic hypergeometric functions obtained by letting $p$ tend to $0$.
\end{abstract}

\maketitle

\section{Introduction}
Elliptic hypergeometric series were introduced by Frenkel and Turaev \cite{FT}. They are a generalization of
the basic and ordinary hypergeometric series. During the last decade several important identities involving 
elliptic hypergeometric series have been found. Typically these identities only involve terminating series, as 
otherwise there would arise very complicated convergence conditions. There also exist elliptic hypergeometric 
integrals, which, when properly specialized, reduce to elliptic hypergeometric series. 

Various identities for these elliptic hypergeometric integrals exist as well. 
These include Spiridonov's elliptic beta integral evaluation \cite{Spir1}, 
and the beta integral transformation formula \cite{Spir2}. The transformation formula was shown \cite{Rainstrafo} to have a Weyl 
group symmetry of type $E_7$, by which we mean that there exists a faithful action of the Weyl group $W(E_7)$ on the parameter 
space of the beta integral which leaves the integral invariant. 

This year, Langer, Schlosser, and Warnaar \cite{LSW} obtained a new elliptic hypergeometric series transformation. 
The main result of this article
is to give an elliptic hypergeometric integral transformation generalizing this series transformation. The Weyl group associated to 
this transformation is of type $F_4$. As far as the author knows this is the first transformation with such a symmetry structure. 
The associated integral can be written as an elliptic beta integral with $16$ parameters, 5 of which are independent.

In a recent article \cite{vdBR} by Rains and the author  many basic hypergeometric identities were obtained from 
the elliptic hypergeometric beta integral evaluation and the transformation with $W(E_7)$ symmetry mentioned above, by
taking a proper limit. In the same vein we consider the limits of this new transformation. In this way we obtain, amongst other things, 
a basic hypergeometric integral with the same $W(F_4)$ symmetry, and an integral expression for a 
very well poised ${}_{14}W_{13}$, both of these results appear to be new. 

All previously known identities mentioned in this introduction have multivariate counterparts. Indeed the proof of the 
series identity from \cite{LSW} is so deeply rooted in multivariate theory that the authors
offered a reward for a proof of the univariate case using just univariate series identities. In contrast, the proof of its univariate integral generalization presented in this article is based
purely on univariate integral identities\footnote{I leave it to the authors of \cite{LSW} to decide if this is sufficient to claim the reward}, and 
we can not prove a multivariate extension. It should be noted that it seems a multivariate integral transformation generalizing the multivariate series identity from \cite{LSW} would not preserve the $W(F_4)$ symmetry.

The organization of this article is as follows: First we have a section on basic notations and definitions. In the next section 
we state and prove the main theorem and show it leads to a $W(F_4)$ symmetry. In the final section we consider the basic hypergeometric
limits.

\subsection{Acknowledgments} I would like to thank Ole Warnaar for pointing out the new elliptic hypergeometric series transformation to me and for our interesting discussions. I would also like to thank Eric Rains for his comments.

\section{Notation}
In this section we give the necessary notation. We use the standard notation from \cite{GR}.
 We assume throughout that $|p|,|q|<1$, which implies all infinite products converge.

We define the $q$-shifted factorials as 
\[
(x;q)_m = \prod_{r=0}^{m-1} (1-xq^r), \qquad (x;q) = (x;q)_\infty = \prod_{r=0}^\infty (1-xq^r).
\]
We use the standard notational abbreviations for $q$-shifted factorials and related functions, writing for example
\[
(a_1,\ldots, a_k;q) = \prod_{r=1}^k (a_r;q), \qquad (a x^{\pm 1};q) = (ax,a/x;q).
\]
Basis hypergeometric series are defined as 
\[
\rphisx{r+1}{r}{a_1,a_2,\ldots,a_{r+1} \\ b_1,b_2,\ldots,b_r}{z}  =
\sum_{k=0}^\infty \frac{(a_1,a_2,\ldots,a_{r+1};q)_k}{(q,b_1,b_2,\ldots,b_r;q)_k} z^k, 
\]
and the very well poised series are given by
\[
\rWs{r+1}{r}(a;b_1,\ldots,b_{r-2};q,z) := 
\rphisx{r+1}{r}{a, \pm q \sqrt{a}, b_1,\ldots,b_{r-2} \\ 
\pm \sqrt{a},aq/b_1,\ldots,aq/b_{r-2}}{z}.
\]
We define the theta functions by
\[
\theta(x;p)= (x,p/x;p). 
\]
The $p,q$-shifted factorials and the elliptic gamma function are defined as
\[
(x;p,q) = \prod_{r,s=0}^{\infty} (1-xp^rq^s), \qquad \Gamma(x) = \Gamma(x;p,q) = \frac{(pq/x;p,q)}{(x;p,q)}.
\]
All gamma functions in this article are elliptic gamma functions, so there should not arise confusion with Euler's gamma function.
We would like to note the following analogs of Legendre's duplication formula, which are used throughout the paper
\begin{align*}
(\pm \sqrt{z}, \pm \sqrt{qz};q)_k & = (z;q)_{2k},  &
(\pm \sqrt{z},\pm \sqrt{pz}, \pm \sqrt{qz}, \pm \sqrt{pqz};p,q) &= (z;p,q), \\
 (\pm \sqrt{z}, \pm \sqrt{qz};q) &= (z;q), &
\Gamma(\pm \sqrt{z},\pm \sqrt{pz}, \pm \sqrt{qz}, \pm \sqrt{pqz}) &= \Gamma(z). 
\end{align*}
Moreover the elliptic gamma function satisfies the difference and reflection relations
\begin{align}\label{eqellgameq}
\Gamma(px) &= \theta(x;q) \Gamma(x), &
\Gamma(x, pq/x) &= 1.
\end{align}

\subsection{Elliptic beta integrals}
We use the definition of the elliptic beta integrals from \cite{vdBR}, and more details can be found there. In particular, the reader should note the prefactor which
ensures that the functions $E^m$ are holomorphic (in $t_r$, $p$, and $q$). 
\begin{definition}
Let $m \in \mathbb{Z}_{\geq 0}$.  Define the set $\mathcal{H}_m =\{ z\in
\mathbb{C}^{2m+6} ~|~ \prod_{i} z_i = (pq)^{m+1} \}/ \sim$, where $\sim$ is
the equivalence relation induced by $z\sim -z$. For parameters $t\in
\mathcal{H}_m$ we define the renormalized elliptic beta integral by
\begin{equation}\label{eqdefem}
E^m(t) =\Biggl(\prod_{0\leq r<s\leq 2m+5} (t_rt_s;p,q)\Biggr)
  \frac{(p;p) (q;q)}{2} \int_{\mathcal{C}} \frac{\prod_{r=0}^{2m+5} \Gamma(t_r z^{\pm 1})}{\Gamma(z^{\pm 2})} \frac{dz}{2\pi i z}.
\end{equation}
where the integration contour $\mathcal{C}$ circles once around the origin
in the positive direction and separates the poles at $z = t_r p^j q^k$
($0\leq r\leq 2m+5$ and $j,k\in \mathbb{Z}_{\geq 0}$) from the poles at
$z=t_r^{-1} p^{-j}q^{-k}$ ($0\leq r\leq 2m+5$ and $j,k \in \mathbb{Z}_{\geq
  0}$).  For parameters $t$ for which such a contour does not exist
(i.e. if $t_rt_s \in p^{\mathbb{Z}_{\leq 0}} q^{\mathbb{Z}_{\leq 0}}$) we
define $E^m$ to be the analytic continuation of the function to these
parameters. 
\end{definition}

One of the most important results about these elliptic beta integrals is the evaluation formula \cite{Spir1} for the $E^0$
\begin{equation}\label{eqeval}
E^0(t) = \prod_{0\leq r<s\leq 5} (pq/t_rt_s;p,q).
\end{equation}
Moreover the integral $E^1$ satisfies a Weyl group of $E_7$ symmetry, which, apart from the permutation symmetry of its eight parameters $t$, just means that it satisfies the equation \cite{Spir2}
\begin{equation}\label{eqe1trafo}
E^1(t) = E^1(t_0v,t_1v,t_2v,t_3v,t_4/v,t_5/v,t_6/v,t_7/v),
\end{equation}
where $v^2=pq/t_0t_1t_2t_3=t_4t_5t_6t_7/pq$.

\subsection{The Weyl group of type $F_4$}
For a thorough introduction to Weyl groups see \cite{Hump}. Roughly speaking; Weyl groups describe finite configurations of mirrors in space, which remain invariant if you look in any one mirror. There are a limited number of configurations like this in a space of given dimension, and $F_4$ refers to one particular case in 4-dimensional space.
\begin{definition}
Let $e_j$ ($j=1,2,3,4$) denote a standard normal basis for $\mathbb{R}^4$.
The root system $R(F_4)$ of type $F_4$ is given by the 48 vectors
\[
R(F_4) = \{\pm e_j ~(1\leq j\leq 4), \frac12 (\pm e_1 \pm e_2 \pm e_3 \pm e_4) , \pm e_j \pm e_k~ (1\leq j<k\leq 4)\}.
\]
For any $\alpha \in R(F_4)$, the reflection $s_{\alpha}$ is given by $s_{\alpha}(\beta) = \beta - 2 \frac{\langle \alpha, \beta\rangle}{\langle \alpha, \alpha \rangle} \alpha$.
The Weyl group $W(F_4)$ of type $F_4$ is the group generated by the reflections $\{ s_{\alpha} ~|~ \alpha \in R(F_4)\}$. A basis for
$F_4$ is given by $\Delta = \{\epsilon_2 -\epsilon_3, \epsilon_1-\epsilon_2, -\epsilon_1, 
\frac12 (\epsilon_1+\epsilon_2+\epsilon_3+\epsilon_4)\}$ (which implies that $W(F_4)$ is generated by the four reflections
$s_{\delta}$ for $\delta \in \Delta$). 
\end{definition}
We recall that the Dynkin diagram of $F_4$ is given by 
\begin{center}
\begin{picture}(110,20)
\put(10,10){\circle*{5}}
\put(40,10){\circle*{5}}
\put(70,10){\circle*{5}}
\put(100,10){\circle*{5}}
\drawline(10,10)(40,10) 
\drawline(70,10)(100,10)
\drawline(40,11)(70,11)
\drawline(40,9)(70,9)
\end{picture}
\end{center}
The Weyl group $W(F_4)$ naturally acts (faithfully) on $\mathbb{C}^4$ by considering the reflections $s_{\alpha}$ as complex reflections.
Given any constant $a$, we can use this standard action of $W(F_4)$ to define a multiplicative
action on $\mathbb{C}^4/ \sim$, where $\sim$ is once again the equivalence relation induced by $z\sim -z$.
\begin{definition}
Given a constant $A$ define the multiplicative action of $W(F_4)$ on $\mathbb{C}^4/ \sim$ via 
\[
w(z) = \exp \circ T_A^{-1} \circ w \circ T_A \circ \log (z).
\]
On the right hand side we use $\log(z_1,z_2,z_3,z_4) = (\log(z_1),\log(z_2),\log(z_3),\log(z_4))$ (and similarly for $\exp$), 
$T_A$ denotes the shift $T_A(z) = (z_1- \frac12 \log(A),z_2- \frac12 \log(A),z_3- \frac12 \log(A),z_4- \frac12 \log(A))$, and
$w$ on the right hand side is the standard action of $W(F_4)$ on $\mathbb{C}^4$.
\end{definition}
Note that the definition is independent of
the choices of logarithm (a different choice at most introduces a common factor $-1$, which explains why we modded out by $\sim$).

More explicitly this action for the basis roots of $W(F_4)$ gives
\begin{align*}
s_{\epsilon_1-\epsilon_2} (z_1,z_2,z_3,z_4) &= (z_2,z_1,z_3,z_4), \\
s_{\epsilon_2-\epsilon_3} (z_1,z_2,z_3,z_4) &= (z_1,z_3,z_2,z_4), \\
s_{-\epsilon_1} (z_1,z_2,z_3,z_4) &= (\frac{A}{z_1},z_2,z_3,z_4), \\
s_{\frac12(\epsilon_1+\epsilon_2+\epsilon_3+\epsilon_4)} (z_1,z_2,z_3,z_4) &=( \frac{Az_1}{\sqrt{z_1z_2z_3z_4}},\frac{Az_2}{\sqrt{z_1z_2z_3z_4}},\frac{Az_3}{\sqrt{z_1z_2z_3z_4}},\frac{Az_4}{\sqrt{z_1z_2z_3z_4}}).
\end{align*}
Observe that the multiplicative action of $W(F_4)$ is faithful as well. The subgroup of $W(F_4)$ generated by the reflections $\{s_{\alpha}~|~ \alpha = \pm e_j, \pm e_j \pm e_k (j\neq k)\}$ is the Weyl group $W(B_4)$ of type $B_4$ and has index $[W(F_4):W(B_4)]=3$. The multiplicative action 
of $W(B_4)$ is by permutations and flips of an arbitrary number of parameters.

\section{The New Elliptic Hypergeometric Identity}
In this section we give the main theorems. In particular we introduce the elliptic hypergeometric
integral of interest and prove the new transformation.

Let us begin by defining the following elliptic hypergeometric function.
\begin{definition}\label{defellint}
The elliptic hypergeometric integral $E(b;t;p,q)$ for $b\in \mathbb{C}$ and $t\in \mathbb{C}^4 / \sim$ where $\sim$ is
still the equivalence relation induced by $z\sim -z$, is given by
\[
E(b;t;p,q) := \frac{E^5(t_1,\frac{pq}{bt_1}, t_2, \frac{pq}{bt_2},t_3,\frac{pq}{bt_3},t_4,\frac{pq}{bt_4},
\pm \sqrt{b}, \pm \sqrt{bq}, \pm \sqrt{bp}, \pm \sqrt{bpq};p,q)}{\prod_{r=1}^4 (bt_r^2, \frac{p^2q^2}{bt_r^2};p,q)}.
\]
\end{definition}
We would like to point out that $E$ has at most simple poles at the zeros of the denominator, that is at
$bt_r^2=p^{-k}q^{-l}$ or $\frac{p^2q^2}{bt_r^2} = p^{-k}q^{-l}$ for some $k,l \in \mathbb{Z}_{\geq 0}$. However, the transformation we
will prove below will show that it cannot have poles at these points. Thus we find that $E$ is holomorphic in 
$(b;t;p,q) \in \mathbb{C} \times \mathbb{C}^4 \times D(0,1)^2$ (where $D(0,1) = \{z \in \mathbb{C}~|~|z|<1\}$ is the open unit disc).
We could divide by some more $p,q$-shifted factorials, while preserving holomorphicity of $E$ (thus removing some excessive zeros). However, finding the largest denominator
such that $E$ is holomorphic seems to be non-trivial, thus we refrain from finding this denominator. This has as only consequence that there will be some superfluous constants on both sides of some of the equations given below; as all equations are equation between holomorphic functions, we can always divide by these factors (to obtain an equation between meromorphic functions).

If we write $E^5$ explicitly as integral, and apply some of the identities for the $p,q$-shifted factorials and the elliptic gamma function, we obtain that, for generic parameters, $E$ is given by
\begin{multline}\label{eqdefE2}
E(b;t;p,q) = 
\prod_{1\leq r<s\leq 4} (t_rt_s, \frac{pqt_r}{bt_s}, \frac{pqt_s}{bt_r}, \frac{p^2q^2}{b^2t_rt_s};p,q) 
\frac{(\frac{pq}{b};p,q)^4 (b^2,pb^2,qb^2,pqb^2;p,q)}{(b,pb,qb,pqb;p,q)}
\\ \times \frac{(p;p) (q;q)}{2}
\int_{\mathcal{C}}  \frac{\Gamma(bz^{\pm 2})}{\Gamma(z^{\pm 2})} \prod_{r=1}^4 \frac{\Gamma(t_rz^{\pm 1})}{\Gamma(bt_r z^{\pm 1})}
\frac{dz}{2\pi i z}.
\end{multline}

We can now give the main theorem
\begin{thm}\label{thmmain}
The following transformation holds for the elliptic hypergeometric integral $E$.
\[
E(b;t_1,t_2,t_3,t_4;p,q) = E(b;t_1v,t_2v, t_3v, t_4v ;p,q),
\]
where $v^2=p^2q^2/b^2t_1t_2t_3t_4$.
\end{thm}
Note that the choice of sign of $v$ does not matter. Also observe that the zeros of the denominator of the elliptic hypergeometric integral on the right hand side are indeed disjoint from those on the left hand side. Therefore $E$ is holomorphic in all its parameters. 

For generic parameters we can write this transformation explicitly (using the $t_r \to pq/bt_r$ symmetry on the right hand side to simplify the formula) as 
\[
\frac{(p;p) (q;q)}{2}
\int_{\mathcal{C}}  \frac{\Gamma(bz^{\pm 2})}{\Gamma(z^{\pm 2})} \prod_{r=1}^4 \frac{\Gamma(t_rz^{\pm 1})}{\Gamma(bt_r z^{\pm 1})}
\frac{dz}{2\pi i z} =
\frac{(p;p) (q;q)}{2}
\int_{\mathcal{C}}  \frac{\Gamma(bz^{\pm 2})}{\Gamma(z^{\pm 2})} \prod_{r=1}^4 \frac{\Gamma(\frac{\sqrt{t_1t_2t_3t_4}}{t_r}z^{\pm 1})}{\Gamma(
\frac{b\sqrt{t_1t_2t_3t_4}}{t_r} z^{\pm 1})}
\frac{dz}{2\pi i z}. 
\]
The transformation reduces to \cite[(4.2)]{LSW} upon setting $t_3\cdot pq/bt_4=q^{-n}$ (which changes the integrals into finite sums of residues). 
Just like its series counterpart we cannot obtain an evaluation formula from this transformation by specializing the integral on one side of the equation to an evaluation. If we do so, either the right hand side becomes an evaluation as well, or the equation we obtain is $0=0$ 
(and if we divide by a term before taking the limit in order to avoid the zero, we have to consider derivatives of elliptic gamma functions).

\begin{proof}
Let us define $s_1 = \sqrt{t_3t_4b^2/pq }$, $s_2=\sqrt{p q/t_3t_4}$, $u_1=\sqrt{ pq t_3/b^2t_4}$ and $u_2 = \sqrt{p q t_4/b^2t_3}$, which implies  $s_1s_2=b$, $s_1u_1=t_3$, $s_1u_2=t_4$ and $s_2u_1=pq/bt_4$ and $s_2u_2= pq/bt_3$ (that is, we choose the roots such that these equations hold). 
Moreover we use the abbreviation $K=(p;p)(q;q)/2$. 
In the following calculation all integrals are over the unit circle. In order for this to be the right contour we impose the conditions 
$|t_r|<1$, $|pq/bt_r|<1$, $|b|<1$, $|s_r|<1$, $|u_r|<1$, $|s_1v|<1$, $|s_2/v|<1$, $|t_rv|<1$, $|pq/bvt_r|<1$. This is an open set of conditions which is satisfied at the points $t_1=t_2=q^{3/4}$, $t_3=t_4=q^{1/2}$, $b=q^{3/4}$, $p=q <1$ (which makes $v=1$), so it satisfied in a non-empty open set. By analytic continuation the final result then also holds (as an identity between meromorphic functions) for all parameters 
$(t,b;p,q)$ in $\mathbb{C}^5 \times D(0,1)^2$. 

We now derive the transformation by a straightforward calculation. First we use the evaluation formula \eqref{eqeval} for an $E^0$ to 
obtain
\begin{align*}
 \int \frac{\Gamma(bz^{\pm 2})}{\Gamma(z^{\pm 2})}  
\prod_{r=1}^4 \frac{\Gamma(t_r z^{\pm 1})}{\Gamma(bt_r z^{\pm 1})}  \frac{dz}{2\pi i z} &=
 \int \frac{\Gamma(bz^{\pm 2})}{\Gamma(z^{\pm 2})}  
\prod_{r=1}^4 \Gamma(t_r z^{\pm 1}, \frac{pq}{bt_r} z^{\pm 1})  \frac{dz}{2\pi i z} \\ &=
\frac{K}{\Gamma(b,b, s_1^2,s_2^2, u_1u_2)} \iint \frac{\Gamma(s_1 z^{\pm 1} y^{\pm 1},s_2 z^{\pm 1} y^{\pm 1}) \Gamma(u_1 y^{\pm 1}, u_2 y^{\pm 1}) }{\Gamma(y^{\pm 2}, z^{\pm 2})} 
\\ & \qquad \qquad \qquad \qquad \qquad \qquad \times 
\Gamma(t_1 z^{\pm 1}, \frac{pq}{bt_1} z^{\pm 1}, t_2 z^{\pm 1}, \frac{pq}{bt_2} z^{\pm 1})
\frac{dy}{2\pi i y} \frac{dz}{2\pi i z}.
\end{align*}
Using the $W(E_7)$ symmetry \eqref{eqe1trafo} in the $z$-integral, where the parameters $t_0,t_1,t_2,t_3$ in \eqref{eqe1trafo} are chosen to be $(s_1y, s_1/y, t_1,t_2)$, and subsequently simplifying using the reflection equation \eqref{eqellgameq} for the elliptic gamma function with $s_1t_1 \cdot s_2pq/bt_1=pq$ (and likewise for $t_2$) shows the original integral equals
\begin{align*}
&=\frac{K\Gamma(s_1^2,t_1t_2,s_2^2, p^2q^2/b^2t_1t_2)}{\Gamma(b,b, s_1^2,s_2^2, u_1u_2)}  \iint \frac{\Gamma(s_1v z^{\pm 1} y^{\pm 1},\frac{s_2}{v} z^{\pm 1} y^{\pm 1})  }{\Gamma(y^{\pm 2}, z^{\pm 2})} \\ 
& \qquad \qquad \qquad \qquad \qquad \qquad \qquad \qquad \times \Gamma(u_1 y^{\pm 1}, u_2 y^{\pm 1},s_1 t_1y^{\pm 1}, s_1 t_2 y^{\pm 1}, \frac{s_2 pq}{bt_1} y^{\pm 1}, \frac{s_2 pq}{bt_2} y^{\pm 1})
\\ & \qquad \qquad \qquad \qquad \qquad \qquad \qquad \qquad \qquad \qquad \qquad \times \Gamma(t_1 v z^{\pm 1}, \frac{pb}{v t_1} z^{\pm 1}, t_2 v z^{\pm 1}, \frac{pq}{bvt_2} z^{\pm 1})   \frac{dy}{2\pi i y} \frac{dz}{2\pi i z} 
\\&= \frac{K\Gamma(t_1t_2, p^2q^2/b^2t_1t_2)}{\Gamma(b,b, u_1u_2)}  \iint \frac{\Gamma(s_1v z^{\pm 1} y^{\pm 1},\frac{s_2}{v} z^{\pm 1} y^{\pm 1}) \Gamma(u_1 y^{\pm 1}, u_2 y^{\pm 1}) }{\Gamma(y^{\pm 2}, z^{\pm 2})} \\ & \qquad \qquad \qquad \qquad \qquad \qquad \qquad \qquad \times \Gamma(t_1 v z^{\pm 1}, \frac{pq}{bv t_1} z^{\pm 1}, t_2 v z^{\pm 1}, \frac{pq}{bvt_2} z^{\pm 1})   \frac{dy}{2\pi i y} \frac{dz}{2\pi i z}.
\end{align*}
Now we can once again use the evaluation formula for an $E^0$ in the $y$-integral, and subsequently simplify again to obtain that this equals
\begin{align*}
 &= \frac{\Gamma(t_1t_2, p^2q^2/b^2t_1t_2, s_1^2 v^2, s_2^2/v^2, b,b,u_1u_2)}{\Gamma(b,b, u_1u_2)} \\ & \qquad \times \int \frac{\Gamma(b z^{\pm 2}) \Gamma( t_3 v z^{\pm 1}, t_4 v z^{\pm 1},\frac{pq}{bvt_4} z^{\pm 1}, \frac{pq}{bvt_3} z^{\pm 1},t_1 v z^{\pm 1}, \frac{pq}{bv t_1} z^{\pm 1}, t_2 v z^{\pm 1}, \frac{pq}{bvt_2} z^{\pm 1}) }{\Gamma(z^{\pm 2})} \frac{dz}{2\pi i z}
\\ &=  \int \frac{\Gamma(b z^{\pm 2})}{\Gamma(z^{\pm 2})}  \prod_{r=1}^4 \Gamma( t_r v z^{\pm 1}, \frac{pq}{bvt_r} z^{\pm 1}) \frac{dz}{2\pi i z}.
\end{align*}
Multiplying the thus obtained integral identity with the correct prefactor (from \eqref{eqdefE2}, note the required prefactor is the same on both sides of the equation) we get the desired transformation.
\end{proof}

We can now easily relate the transformation for $E$ with the multiplicative $W(F_4)$ action.
\begin{cor}\label{corwf4}
Let $W(F_4)$ act on $\mathbb{C}^4/ \sim$ using the multiplicative action with parameter $A=pq/b$. Then we have
\[
E(b;t;p,q) = E(b;w(t);p,q)
\]
for any $w\in W(F_4)$.
\end{cor}
\begin{proof}
It suffices to show that the equation holds for a set of generators of $W(F_4)$, thus we need to show it for
$s_{\delta}$ ($\delta \in \Delta$). The action of $s_{\frac12(\epsilon_1+\epsilon_2+\epsilon_3+\epsilon_4)} $ is exactly the transformation of
Theorem \ref{thmmain}. Moreover the actions of $s_{-\epsilon_1}$, $s_{\epsilon_1-\epsilon_2}$ and $s_{\epsilon_2-\epsilon_3}$ just permute 
the arguments of $E^5$ in Definition \ref{defellint}, so they clearly preserve $E$. 
\end{proof}

\section{Basic Hypergeometric Limits}
In this section we consider some limits as $p\to 0$. In particular this will give us some new basic hypergeometric identities. 

Before we start taking the limit as $p\to 0$ we need to specify how our other parameters depend on $p$ (that is, how they behave as $p\to 0$). 
In all cases we will assume $q$ is fixed and the other parameters depend on $p$ via some power. That is we consider 
\[
E(b p^{\beta}; t_1 p^{\tau_1}, t_2 p^{\tau_2}, t_3p^{\tau_3}, t_4p^{\tau_4};p,q).
\]
Note that the $W(F_4)$ symmetry of $E$ becomes a $W(F_4)$ action on $(\beta,\tau)$, which leaves $\beta$ invariant and is
the standard action on $\mathbb{R}^4$ with reflections on the $\tau$-variables. 

There are some very simple limits. Indeed if we just use the elementary limits
\[
\lim_{p\to 0} (p^{\alpha} z;p,q) = \begin{cases}
(z;q) & \alpha=0, \\
1 & 0<\alpha,
\end{cases}
\]
and hence,
\[
\lim_{p\to 0} \Gamma(p^{\alpha} z;p,q) = \begin{cases}
\frac{1}{(z;q)} & \alpha=0, \\
1 & 0<\alpha<1, \\
(q/z;q) & \alpha=1,
\end{cases}
\]
we can take limits directly in \eqref{eqdefE2} if $0\leq \beta \leq 1$ and $0\leq \tau_r \leq 1-\beta$ ($r=1,2,3,4$).
In this article we therefore consider the $W(F_4)$ orbit of this set, which is the polytope with bounding equations 
$0\leq \beta\leq 1$, $0\leq \tau_r+\tau_s \leq 2-2\beta$ (for $r\neq s$), and $\tau_r-\tau_s \leq 1-\beta$ (for $r\neq s$). 
This polytope is a pyramid with the 24-cell as base (in the plane $\beta=0$) and the point at $\beta=1$ and $\tau=(0,0,0,0)$ as its apex.

As the transformation leaves $b$ and thus $\beta$ invariant, it is convenient to order the different limits by the different values of $\beta$. 
\subsection{$\beta=1$}
The only point of interest here is $\tau_r=0$ for all $r$. 
Let us first define (in similar vein to the definition of the elliptic beta integral) the following basic hypergeometric integral.
\begin{definition} For $b\in \mathbb{C}$, $t \in \mathbb{C}^4/\sim$ and $|q|<1$ we set
\[
B_1(b;t;q) := 
\prod_{1\leq r<s\leq 4} (t_rt_s, \frac{qt_r}{bt_s}, \frac{qt_s}{bt_r}, \frac{q^2}{b^2t_rt_s};q) 
(\frac{q}{b};q)^4 
\frac{ (q;q)}{2}
\int_{\mathcal{C}}  \frac{\Gamma(\frac{q}{b}z^{\pm 2},z^{\pm 2};q)}{ \prod_{r=1}^4 \Gamma(t_rz^{\pm 1},\frac{q}{bt_r} z^{\pm 1};q)}
\frac{dz}{2\pi i z},
\]
where $\mathcal{C}$ is a deformation of the unit circle traversed in positive direction, which separates the poles at
$z=t_r q^k$ and $z=\frac{1}{bt_r} q^{k+1}$ ($1\leq r\leq 4$ and $k\in \mathbb{Z}_{\geq 0}$) from their reciprocals. If such a contour does not exist, then $B_1$  is equal to the analytic continuation of the integral on the right hand side (in $b$ and $t_r$) to the desired parameters.
\end{definition}
This definition was inspired by the following result
\begin{prop}
We have
\[
\lim_{p\to 0} E(pb;t;p,q) = 
B_1(b;t;q).
\]
\end{prop}
\begin{proof}
This is just a special case of \cite[Proposition 1]{vdBR}. It can also be obtained by taking the limit in \eqref{eqdefE2} and interchanging limit and integral.
\end{proof}
Note that the basic hypergeometric integral on the right hand side is now seen to satisfy a $W(F_4)$ symmetry. 
\begin{cor}
Let $W(F_4)$ act on $\mathbb{C}^4 /\sim$ using the multiplicative action with parameter $A=q/b$. Then we have
\[
B_1(b;t;q) = B_1(b;w(t);q)
\]
for any $w\in W(F_4)$.
\end{cor}
\begin{proof}
Just take the limit $p\to 0$ in Corollary \ref{corwf4} with $b$ replaced by $pb$.
\end{proof}

\subsection{$0<\beta<1$}
The limits in this region are quite boring. Indeed, if we restrict to the polytope $0\leq \tau_r \leq 1-\beta$, we see that the limit
in \eqref{eqdefE2} gives (a special case of) the Askey-Wilson integral \cite[(6.1.1)]{GR}, which satisfies an evaluation. Unsurprisingly
this means that in fact we just get an evaluation in the entire polytope.

\begin{prop}
Let $0< \beta<1$ and $0\leq \tau_r+\tau_s \leq 2-2\beta$ (for $r\neq s$) and 
$\tau_r-\tau_s \leq 1-\beta$ ($r\neq s$). We obtain the following limits
\begin{itemize}
\item If all four $\tau_r$ are either $0$ or $1-\beta$ we have
\[
\lim_{p\to 0} E(bp^{\beta}; t_r p^{\tau_r};p,q) = ( \prod_{r: \tau_r=0} t_r \prod_{r: \tau_r=1-\beta} \frac{q}{bt_r};q);
\]
\item If $\tau_1=\tau_2=\tau_3=(1-\beta)/2$ and $\tau_4=(3-\beta)/2$ we have
\[
\lim_{p\to 0} E(bp^{\beta}; t_r p^{\tau_r};p,q) = ( \frac{q^3}{b^3 t_4^2};q);
\]
\item If $\tau_1=\tau_2=\tau_3=(1-\beta)/2$ and $\tau_4=(\beta-1)/2$ we have
\[
\lim_{p\to 0} E(bp^{\beta}; t_r p^{\tau_r};p,q) = ( \frac{q t_4^2}{b};q);
\]
\item Otherwise (that is, in those cases in which none of the permutations of the $\tau_r$ are covered above) we have 
\[
\lim_{p\to 0} E(bp^{\beta}; t_r p^{\tau_r};p,q) = 1.
\]
\end{itemize}
\end{prop}
\begin{proof}
Within the polytope $0\leq \tau_r\leq 1-\beta$, this is a special case of \cite[Lemma 2]{vdBR}. In that case we can also obtain the result by exchanging limit and integral in \eqref{eqdefE2} and using the Askey-Wilson integral evaluation to evaluate the limit. 
The simplest way to obtain the other limits is to use the $W(F_4)$ symmetry on the left hand side to map any case to the polytope $0\leq \tau_r\leq 1-\beta$ and then apply the known limit. 

A different method (which avoids using the $W(F_4)$ symmetry) in the case $\tau_4<0$ would be to note that the contour in the original integral 
has to include (in its interior) the poles at $z=t_rp^{\tau_r} q^k$ ($k\in \mathbb{Z}_{\geq 0}$), which (for fixed $k$) move to infinity as $p\to 0$. All other poles however will either be on the right side of the unit circle as $p\to 0$, or remain fixed as $p \to 0$. In order to take the limit it is thus opportune to rewrite the integral as the sum of residues at $z=t_rp^{\tau_r} q^k$ (and their reciprocals) and an integral with integration contour close to the unit circle (that is, as $p\to 0$ we pick up more and more poles to keep the contour roughly constant). Subsequently we can bound the integral and poles in such a way that it can be shown only the residue at $z=t_rp^{\tau_r}$ (and 
symmetrically $z=t_r^{-1} p^{-\tau_r}$) contributes to the limit. The limit of this residue is then easily calculated.
Similarly one can treat the case $\tau_4>1-\beta$. The explicit calculations are quite tedious, and considering we already have a simple argument proving the proposition we omit the details.
\end{proof}

\subsection{$\beta=0$}
In this case we need several expressions to describe all limits. Let us begin with the simplest case.
\begin{prop}\label{propb01}
If $0\leq \tau_r \leq 1$ for $1\leq r\leq 4$ we have
\begin{multline*}
\lim_{p\to 0} E(b; t_rp^{\tau_r};p,q) = 
\prod_{r<s:\tau_r=\tau_s=0} (t_rt_s;q)
\prod_{r,s: \tau_r=0, \tau_s=1} (\frac{qt_r}{bt_s};q)
\prod_{r<s: \tau_r=\tau_s=1} (\frac{q^2}{b^2t_rt_s};q)
\frac{(b^2,qb^2;q)}{(b,qb;q)}
\\ \times \frac{ (q;q)}{2}
\int_{\mathcal{C}}  \frac{(z^{\pm 2};q)}{(bz^{\pm 2};q)} 
\prod_{r:\tau_r=0} \frac{(bt_rz^{\pm 1};q)}{(t_rz^{\pm 1};q)}
\prod_{r:\tau_r=1} \frac{(\frac{q}{t_r} z^{\pm 1};q)}{(\frac{q}{bt_r} z^{\pm 1};q)}
\frac{dz}{2\pi i z}.
\end{multline*}
where $\mathcal{C}$ is a deformation of the unit circle traversed in positive direction, which separates the zeros 
of terms of the form $(xz;q)$ (for $x=t_r$ ($\tau_r=0$) and $x=q/bt_r$ ($\tau_r=1$)) and $(bz^2;q)$ from their reciprocals. If such a contour does not exist, then the limit is equal to the analytic continuation of the integral on the right hand side (in $b$ and $t_r$) to the desired parameters.
\end{prop}
\begin{proof}
This is again a special case of \cite[Proposition 1]{vdBR}; once again we can also obtain it by taking the limit in \eqref{eqdefE2} and interchanging limit and integral.
\end{proof}

For the rest of the limits we obtain
\begin{prop}\label{propb02}
If $\tau_1< 0$ and $-\tau_1\leq  \tau_r\leq 1+\tau_1$ we have
\begin{multline}\label{eqlim2}
\lim_{p\to 0} E(b;t_r p^{\tau_r};p,q) =
\prod_{r\geq 2: \tau_r=-\tau_1} (t_rt_1;q)
\prod_{r\geq 2: \tau_r=1+\tau_1} (\frac{qt_1}{bt_r};q)
\frac{ (\frac{q}{b},qb^2;q)}{(\frac{q}{b^2},qb;q)}
(q;q) 
\\  \qquad \qquad \times \int_{\mathcal{C}}  
\left(\frac{(1-z^{-2}) (\frac{q}{bz^2},\frac{bt_1}{z};q)}{(\frac{b}{z^2},\frac{t_1}{z};q)}  \right)^{1_{\{\tau_1=-1/2\}}}
\frac{\theta(b^2t_1z;q)}{(t_1z, \frac{q}{bt_1z};q)} 
\prod_{r\geq 2:\tau_r=-\tau_1} \frac{(\frac{bt_r}{z};q)}{(\frac{t_r}{z};q)}
\prod_{r\geq 2:\tau_r= 1+\tau_1} \frac{(\frac{q}{t_rz};q)}{(\frac{q}{bt_rz};q)}
 \frac{dz}{2\pi i z}.
\end{multline}
where $\mathcal{C}$ is a deformation of the unit circle traversed in positive direction, which separates the zeros 
of $(t_1z;q)$ from those of $(b/z^2,t_1/z;q)$ (only if $\tau_1=-1/2$), $(q/bt_1z;q)$, $(t_r/z;q)$ ($\tau_r=-\tau_1$), and $(q/bt_rz;q)$ ($\tau_r=1+\tau_1$). If such a contour does not exist, then the limit is equal to the analytic continuation of the integral on the right hand side (in $b$ and $t_r$) to the desired parameters.
\end{prop}
\begin{proof}
The idea of the proof is similar to the proof of \cite[Proposition 2]{vdBR}. Indeed we note that 
\[
\frac{\theta(bwz,\frac{w}{bz}, \frac{b}{z^2};q)}{\theta(w z^{\pm 1},\frac{1}{z^2};q)} + (z \to 1/z) = 
\frac{\theta(b^2;q)}{\theta(b;q)},
\]
which is a reformulation of Riemann's addition formula for theta functions. Noting that the integrand in \eqref{eqdefE2} is even, we can use this equation to obtain
\begin{multline*}
E(b;t;p,q) = 
\prod_{1\leq r<s\leq 4} (t_rt_s, \frac{pqt_r}{bt_s}, \frac{pqt_s}{bt_r}, \frac{p^2q^2}{b^2t_rt_s};p,q) 
\frac{(\frac{pq}{b};p,q)^4 (b^2,pb^2,qb^2,pqb^2;p,q)}{(b,pb,qb,pqb;p,q)}
\\ \times (p;p) (q;q) \frac{\theta(b;q)}{\theta(b^2;q)}
\int_{\mathcal{C}}  \frac{\Gamma(bz^{\pm 2})}{\Gamma(z^{\pm 2})} \prod_{r=1}^4 \frac{\Gamma(t_rz^{\pm 1})}{\Gamma(bt_r z^{\pm 1})}
\frac{\theta(bwz,\frac{w}{bz}, \frac{b}{z^2};q)}{\theta(w z^{\pm 1},\frac{1}{z^2};q)}
\frac{dz}{2\pi i z}.
\end{multline*}
Specializing $w=bt_1$ and using the difference equation \eqref{eqellgameq} for the elliptic gamma function we get
\begin{multline*}
E(b;t;p,q) = 
\prod_{1\leq r<s\leq 4} (t_rt_s, \frac{pqt_r}{bt_s}, \frac{pqt_s}{bt_r}, \frac{p^2q^2}{b^2t_rt_s};p,q) 
\frac{(\frac{pq}{b};p,q)^4 (\frac{q}{b},pb^2,qb^2,pqb^2;p,q)}{(\frac{q}{b^2},pb,qb,pqb;p,q)}
\\ \times (p;p) (q;q) 
\int_{\mathcal{C}}  \frac{\Gamma(bz^2,pb/z^2)}{\Gamma(z^2, p/z^2 )} \frac{\Gamma(pt_1/z, t_1z)}{\Gamma(pbt_1 z^{\pm 1})} \theta(b^2t_1z;q)\prod_{r=2}^4 \frac{\Gamma(t_rz^{\pm 1})}{\Gamma(bt_r z^{\pm 1})}
 \frac{dz}{2\pi i z}.
\end{multline*}
If we plug in $t_r \to t_r p^{\tau_r}$ and subsequently shift the integration variable $z\to zp^{-\tau_1}$ we obtain
\begin{align*}
E(b;t_r p^{\tau_r} ;p,q) &= 
\prod_{1\leq r<s\leq 4} (t_rt_sp^{\tau_r+\tau_s}, \frac{qt_r}{bt_s} p^{1+\tau_r-\tau_s}, \frac{qt_s}{bt_r}p^{1+\tau_s-\tau_r}, \frac{q^2}{b^2t_rt_s}p^{2-\tau_r-\tau_s};p,q) \\ & \qquad \qquad \times 
\frac{(\frac{pq}{b};p,q)^4 (\frac{q}{b},pb^2,qb^2,pqb^2;p,q)}{(\frac{q}{b^2},pb,qb,pqb;p,q)}
(p;p) (q;q) 
\\ & \qquad \qquad \times \int_{\mathcal{C}}  \frac{\Gamma(bz^2p^{-2\tau_1},bp^{1+2\tau_1}/z^2)}{\Gamma(z^2p^{-2\tau_1}, p^{1+2\tau_1}/z^2 )} \frac{\Gamma(t_1p^{1+2\tau_1}/z, t_1z)}{\Gamma(bt_1pz, bt_1p^{1+2\tau_1}/z)} \theta(b^2t_1z;q)
\\ & \qquad \qquad \qquad \qquad \times 
\prod_{r=2}^4 \frac{\Gamma(t_rp^{\tau_r-\tau_1}z, t_rp^{\tau_r+\tau_1}/z )}{\Gamma(bt_r p^{\tau_r-\tau_1}z, b t_r p^{\tau_r+\tau_1}/z )}
 \frac{dz}{2\pi i z}.
 \end{align*}
In the final integral we can now directly take the limit $p\to 0$ of the integrand (as $-1/2\leq \tau_1<0$), moreover the conditions on the contour are such that the contour can remain fixed (to a contour which works for the right hand side of \eqref{eqlim2}) when letting $p$ tend to zero, so the result just follows from plugging in $p=0$.

If no desired contour exists for the right hand side of \eqref{eqlim2} we can use the same argument as in \cite[Proposition 1]{vdBR} to show that the result actually holds as an identity between holomorphic functions, thus in particular also in this case.
\end{proof}

Finally we have the following limit
\begin{cor}\label{corb0}
If $\tau_1> 1$ and $\tau_1-1\leq  \tau_r\leq 2-\tau_1$ we have
\begin{multline}\label{eqlim3}
\lim_{p\to 0} E(b;t_r p^{\tau_r};p,q) =
\prod_{r\geq 2: \tau_r=\tau_1-1} (\frac{qt_r}{bt_1};q)
\prod_{r\geq 2: \tau_r=2-\tau_1} (\frac{q^2}{b^2t_1t_r};q)
\frac{ (\frac{q}{b},qb^2;q)}{(\frac{q}{b^2},qb;q)}
(q;q) 
\\  \qquad \qquad \times \int_{\mathcal{C}}  
\left(\frac{(1-z^{-2}) (\frac{q}{bz^2},\frac{q}{t_1z};q)}{(\frac{b}{z^2},\frac{q}{bt_1z};q)}  \right)^{1_{\{\tau_1=3/2\}}}
\frac{\theta(\frac{qbz}{t_1};q)}{(\frac{q}{bt_1}z, \frac{t_1}{z};q)} 
\prod_{r\geq 2:\tau_r=\tau_1-1} \frac{(\frac{bt_r}{z};q)}{(\frac{t_r}{z};q)}
\prod_{r\geq 2:\tau_r= 2-\tau_1} \frac{(\frac{q}{t_rz};q)}{(\frac{q}{bt_rz};q)}
 \frac{dz}{2\pi i z}.
\end{multline}
where $\mathcal{C}$ is a deformation of the unit circle traversed in positive direction, which separates the zeros 
of $(t_1/z;q)$ from the other zeros of the denominator. If such a contour does not exist, then the limit is equal to the analytic continuation of the integral on the right hand side (in $b$ and $t_r$) to the desired parameters.
\end{cor}
\begin{proof}
Use the $t_1p^{\tau_1} \to \frac{pq}{bt_1p^{\tau_1}}$ symmetry of $E$ and then apply the previous proposition. One could also use a very similar argument to the proof of that proposition.
\end{proof}
The integrals in Proposition \ref{propb02} and Corollary \ref{corb0} have a series representation. It can be directly seen that the corresponding
limits equal those series, by writing the elliptic hypergeometric integral as a sum of residues plus an integral with a different contour, and then taking the limit. However showing convergence in that method is much more complicated than first showing the limit is an integral and then 
calculating the relevant series representation from that integral as in the next proposition.
\begin{prop}
If $|b|<1$ we have 
\begin{align*}
\prod_{r=2}^4 &(t_rt_1, \frac{qt_1}{bt_r};q)
\frac{ (\frac{q}{b},qb^2;q)}{(\frac{q}{b^2},qb;q)}
(q;q) 
 \int_{\mathcal{C}}  
\frac{(1-z^{-2}) (\frac{q}{bz^2},\frac{bt_1}{z};q)}{(\frac{b}{z^2},\frac{t_1}{z};q)} 
\frac{\theta(b^2t_1z;q)}{(t_1z, \frac{q}{bt_1z};q)} 
\prod_{r=2}^4  \frac{(\frac{bt_r}{z}, \frac{q}{t_rz};q)}{(\frac{t_r}{z}, \frac{q}{bt_rz};q)}
 \frac{dz}{2\pi i z} \\ &=
\prod_{r=2}^4 (bt_1t_r, \frac{qt_1}{t_r};q)  
\frac{(qb^2,b^2, \frac{qt_1^2}{b};q)}{(qb, qt_1^2;q)} {}_{14}W_{13}\left(t_1^2; t_1t_2,\frac{qt_1}{bt_2}, t_1t_3, \frac{qt_1}{bt_3}, t_1t_4, \frac{qt_1}{bt_4}, \frac{q}{b}, \pm \sqrt{b} t_1, \pm \sqrt{bq} t_1;q,b^2\right).
\end{align*}
Moreover we have for $k\in \mathbb{Z}_{\geq 0}$ and $|b|<1$
\begin{align*}
\prod_{r=1}^k & (t_1u_r;q) \frac{(\frac{q}{b}, qb^2, q;q)}{(\frac{q}{b^2}, qb;q)}
\int_{\mathcal{C}} \frac{\theta(b^2t_1z;q)}{(t_1z,\frac{q}{bt_1z};q)} \prod_{r=1}^k \frac{(\frac{bu_r}{z};q)}{(\frac{u_r}{z};q)} \frac{dz}{2\pi i z} \\ &= \prod_{r=1}^k (bu_rt_1;q) \frac{(qb^2,b^2;q)}{(qb;q)}
\rphisx{k+1}{k}{t_1u_1, \ldots, t_1u_k, \frac{q}{b} \\ bu_1t_1, \ldots, bu_kt_1}{b^2}.
\end{align*}

\end{prop}
\begin{proof}
Both identities follow from residue calculus. Indeed if we let the contour go to infinity, we have to pick up the poles at $z=t_1^{-1}q^{-k}$. 
The sum of the values of these residues give the series on the right hand sides. The condition that $|b|<1$ ensures that the series converges and that the value of the integrand at infinity goes to zero (exponentially, so the total value of the integral converges to zero as well).
\end{proof}

Together Propositions \ref{propb01}, \ref{propb02} and Corollary \ref{corb0} together with their analogues for permutations of the $\tau_r$ cover the entire 24-cell we are interested in.
We end this section by tallying the different equations we obtain for basic hypergeometric functions by taking the limit in the 
$W(F_4)$ transformations of the elliptic hypergeometric integral. It is known that the 24-cell has a $W(F_4)$ symmetry. Up 
to this $W(F_4)$ symmetry, it has exactly one vertex, one edge, one 2-dimensional face (a triangle), one 3-dimensional face (an octahedron), and, of course, one interior.

We would like to observe here that all representations of the limits we have are invariant under the $t_r\to pq/bt_r$ symmetry on the left hand side, that is, doing this symmetry and then taking the limit gives exactly the same expression as taking the limit immediately. Similarly the permutation symmetry of the elliptic hypergeometric integral becomes the same permutation symmetry of the limiting integral/series. In particular we may hope to find at most $3=[W(F_4):W(B_4)]$ different expressions (either integral or series representations) for each limit.

\subsubsection{The vertex}
The different $\tau$-vectors we consider on this level are in the orbit of $(0,0,0,0)$. That is they are $(0,0,0,0)$, $(0,0,0,1)$, $(0,0,1,1)$, $(0,1,1,1)$, $(1,1,1,1)$, $(-\frac12, \frac12,\frac12,\frac12)$, $(\frac12,\frac12,\frac12,\frac32)$ and permutations of these.
If we consider (for $t\in \mathbb{C}^4$) the integral corresponding to $\tau=(0,0,0,0)$ 
\[
B_2(b;t;q) := \prod_{1\leq r<s\leq 4} (t_rt_s;q) \frac{(b^2,qb^2,q;q)}{2(b,qb;q)}
\int_{\mathcal{C}} \frac{(z^{\pm 2};q)}{(bz^{\pm 2};q)} \prod_{r=1}^4 \frac{(bt_rz^{\pm 1};q)}{(t_rz^{\pm 1};q)} \frac{dz}{2\pi i z},
\]
we observe that $B_2$ is permutation symmetric in the $t_r$ (which is visible from its definition). Also note that we can write all the limits
corresponding to all-integer $\tau$-vectors in terms of $B_2$, however doing that trivializes the $t_r\to pq/bt_r$ reflections, that is, they don't change the parameters of the $B_2$. As a non-trivial symmetry we do get 
\begin{align*}
B_2(b;t_1,t_2,t_3,t_4;q) &= B_2(b;\frac{T}{t_1},\frac{T}{t_2}, \frac{T}{t_3}, \frac{T}{t_4};q),
\end{align*}
where $T=\sqrt{t_1t_2t_3t_4}$ (the choice of root does not matter, as we can set $z\to -z$ in the definition of $B_2$ to multiply all the 
$t$ parameters by $-1$). Moreover we can express it in terms of a ${}_{14}W_{13}$ as follows
\[ 
B_2(b;t;q) = \prod_{1\leq r<s\leq 4} (bt_rt_s;q) \frac{(qb^2, b^2, T^2;q)}{(qb, b T^2;q)}
{}_{14}W_{13}\left( \frac{bT^2}{q};t_1t_2,t_1t_3,t_1t_4,t_2t_3,t_2t_4,t_3t_4,\frac{q}{b}, \pm \frac{bT}{\sqrt{q}}, \pm bT;q, b^2\right).
\] 
\begin{proof}
Indeed, 
\begin{align*}
\lim_{p\to 0} E(b;t_1,t_2,t_3,t_4;p,q) &= \lim_{p\to 0} E(b;\frac{pqt_1}{bT}, \frac{pqt_2}{bT}, \frac{pqt_3}{bT}, \frac{pqt_4}{bT};p,q) \\&
=\lim_{p\to 0} E(b;\sqrt{\frac{bT^2}{pq}}, \sqrt{\frac{pqt_3t_4}{bt_1t_2}}, \sqrt{\frac{pqt_2t_4}{bt_1t_3}}, \sqrt{\frac{pqt_2t_3}{bt_1t_4}};p,q)
\end{align*}
where inside the limits we just used the $W(F_4)$ symmetry of $E$.
\end{proof}
Note that in the integral representation of the ${}_{14}W_{13}$ the series is clearly symmetric under $t_r \to T/t_r$. We don't get any non-trivial symmetries of the ${}_{14}W_{13}$. However, if we were to specialize for example $t_3t_4=q^{-n}$ the integral expression for $B_2(b;t;q)$ would reduce to 
a sum of residues, which forms a terminating ${}_{14}W_{13}$. So in that case we do get a non-trivial transformation for 
terminating ${}_{14}W_{13}$'s. In particular this gives us the $p\to 0$ limit of \cite[(4.2)]{LSW}

\subsubsection{The edge}
We consider the orbit of $\tau=(0,0,0,x)$ for some $0<x<1$. We only get the identity
\begin{multline*}
\prod_{1\leq r<s\leq 3} (t_rt_s;q) \frac{(b^2,qb^2,q;q)}{2(b,qb;q)} \int_{\mathcal{C}} \frac{(z^{\pm 2};q)}{(bz^{\pm 2};q)} \prod_{r=1}^3 \frac{(bt_r z^{\pm 1};q)}{(t_rz^{\pm 1};q)} \frac{dz}{2\pi i z} \\ = 
\prod_{1\leq r<s\leq 3} (bt_rt_s;q) \frac{(qb^2,b^2;q)}{(qb;q)} \rphisx{4}{3}{t_1t_2,t_1t_3,t_2t_3,\frac{q}{b} \\ bt_1t_2,bt_1t_3,bt_2t_3}{b^2}.
\end{multline*}
\begin{proof}
Indeed, for $0<x<1$ we have
\[
\lim_{p\to 0} E(b;t_1,t_2,t_3,p^xu;p,q) = \lim_{p\to 0} E(b;\sqrt{\frac{ut_2t_3}{t_1}}p^{x/2},\sqrt{\frac{ut_1t_3}{t_2}}p^{x/2},\sqrt{\frac{ut_1t_2}{t_3}}p^{x/2},
\sqrt{\frac{t_1t_2t_3}{u}}p^{-x/2};p,q).
\]
\end{proof}
The third coset of $W(B_4)$ in $W(F_4)$ gives the same $\rphis{4}{3}$ representation.

\subsubsection{The triangle}
A prototypical 2-dimensional face has vertices $(0,0,0,0)$, $(0,0,0,1)$, and $(-\frac12, \frac12,\frac12,\frac12)$. 
 Inspection learns that the limits we have give the same expression (a $\rphis{3}{2}$) for the function on this face for all points in 
 the $W(F_4)$ orbit, thus we do not get any equations on this level.

\subsubsection{The octahedron}
The 3-dimensional faces are octahedra, with as typical example the octahedron with vertices $(0,0,0,0)$, $(0,0,0,1)$, $(0,0,1,0)$, $(0,0,1,1)$, $(-\frac12,\frac12,\frac12,\frac12)$, $(\frac12,-\frac12,\frac12,\frac12)$. Almost all limits in this octahedron will give a ${}_2\phi_1$, except the ones on the square with vertices $(0,0,0,0)$, $(0,0,0,1)$, $(0,0,1,0)$, and $(0,0,1,1)$, which gives an integral. Using the $W(F_4)$ symmetry we move to other octahedra which again almost all give the same ${}_2\phi_1$, except in one square, in which we get an integral. The squares, however, are not mapped to each other. So we obtain the following integral representation (while the right hand side looks like the integral representation \cite[(6.4.11)]{GR} for a sum of two very well poised $\rphis{10}{9}$'s, it is not as it does not satisfy the appropriate balancing condition).
\[
\frac{(btv,qb^2,b^2;q)}{(qb;q)} \rphisx{2}{1}{tv,q/b \\ btv}{b^2} = \frac{(tv,b^2,qb^2,q;q)}{2(b,qb;q)}
\int_{\mathcal{C}} \frac{(z^{\pm 2}, btz^{\pm 1}, bvz^{\pm 1};q)}{(bz^{\pm 2},tz^{\pm 1}, vz^{\pm 1};q)} \frac{dz}{2\pi i z}.
\]
\begin{proof}
We have for $0<x<y<1$ 
\begin{multline*}
\lim_{p\to 0} E(b; t,v,u_1 p^{x}, u_2 p^{y};p,q) \\ = \lim_{p\to 0} 
E(b; \sqrt{\frac{u_1u_2v}{t}}p^{(x+y)/2}, \sqrt{\frac{tu_1u_2}{v}}p^{(x+y)/2}, \sqrt{\frac{tvu_2}{u_1}}p^{(y-x)/2}, \sqrt{\frac{tvu_1}{u_2}}p^{(x-y)/2};p,q).
\end{multline*}
\end{proof}

\subsubsection{The interior}
As with the octahedra, different parts of the interior give different expressions (integrals or ${}_1\phi_0$) for the same function, and the $W(F_4)$ symmetry shows they are identical anyway. Two of the expressions are known to have an evaluation, that is, the ${}_1\phi_0$ (obtained from 
$\tau_1<0$ and $-\tau_1<\tau_r<1+\tau_1$ for $r=2,3,4$), and
the Askey-Wilson integral evaluation (obtained from $0<\tau_r<1$ for $r=1,2,3,4$). Using either evaluation shows that the value of the limit in the interior can be expressed as $(qb^2;q)$. For the final expression we find the identity
\[
(qb^2;q) = \frac{(b^2,qb^2,q;q)}{2(b,qb;q)} \int_{\mathcal{C}} \frac{(z^{\pm 2}, btz^{\pm 1};q)}{(bz^{\pm 2}, tz^{\pm 1};q)} \frac{dz}{2\pi i z}.
\]
\begin{proof}
Indeed, for $0<x<2/3$ we have
\[
\lim_{p\to 0} E(b;t,u_1p^{x}, u_2p^{x}, u_3p^{x};p,q) = \lim_{p\to 0}
E(b; \sqrt{\frac{u_1u_2u_3}{t}}p^{3x/2}, \sqrt{\frac{tu_2u_3}{u_1}}p^{x/2}, \sqrt{\frac{tu_1u_3}{u_2}}p^{x/2}, \sqrt{\frac{tu_1u_2}{u_3}}p^{x/2};p,q),
\]
and the right hand side can be evaluated using the Askey-Wilson integral evaluation \cite[(6.1.1)]{GR}. For $2/3<x<1$ the right hand side would become a ${}_1\phi_0$ and we could have used its evaluation instead as well.
\end{proof}
We would like to make some final remarks on this formula. First of all, while the equation looks like the Nassrallah-Rahman integral evaluation, the integral in question does not satisfy the right balancing condition. We can, however, express the integral in terms of a very well poised ${}_8W_7$ and
thus obtain an evaluation for a ${}_8W_7$ with properly specialized values, that is
\[
\frac{(qb^2,b^2t^2;q)}{(b^3t^2,qb;q)} =  {}_8W_7(\frac{b^2t^2}{q}; \pm t\sqrt{b}, \pm t \sqrt{\frac{b}{q}} , b;q,qb)
=\sum_{k=0}^\infty 
\frac{(1-b^2t^2q^{2k-1})}{(1-b^2t^2q^{-1})}
\frac{(b^2t^2/q, b;q)_k
(bt^2/q;q)_{2k}
}{(q, bt^2;q)_k (b^3t^2;q)_{2k}} (bq)^k.
\]
This equation is the special case of \cite[(3.5.3)]{GR} with $b=\sqrt{q}$ (in that equation), where we have to use the $q$-Gauss summation to evaluate
the ${}_2\phi_1$ appearing there.


\begin{thebibliography}{99}
\bibitem{vdBR} F.J. van de Bult and E.M. Rains, \textit{Basic hypergeometric functions as limits of elliptic hypergeometric functions}, SIGMA \textbf{5} (2009), 059, 31 pp, arxiv:0902.0621.

\bibitem{FT} I.B. Frenkel and V.G. Turaev, \textit{Elliptic solutions of the Yang-Baxter equation and modular hypergeometric functions}, The Arnold-Gelfand Mathematical Seminars, Birkh\"auser Boston, Boston, MA, 1997, 171--204.

\bibitem{GR} G. Gasper and M. Rahman, {Basic Hypergeometric Series}, Encyclopedia of Mathematics and its Applications, Vol. 96, 2nd ed., Cambridge University Press, Cambridge, 2004.

\bibitem{Hump} J.E. Humphreys, \textit{Reflection Groups and Coxeter Groups}, Cambridge Studies in Advanced Mathematics \textbf{29}, Cambridge University Press, Cambridge, 1990.

\bibitem{LSW} R. Langer, M.J. Schlosser and S.O. Warnaar, \textit{Theta Functions, Elliptic Hypergeometric Series, and Kawanaka's Macdonald Polynomial Conjecture}, SIGMA \textbf{5} (2009), 055, 20 pp, arxiv:0905.4033. 

\bibitem{Rainstrafo} E.M. Rains, \textit{Transformations of elliptic hypergeometric integrals}, Ann. Math. (2009), to appear.

\bibitem{Spir1} V.P. Spiridonov, \textit{On the elliptic beta function}, Russian Math. Surveys \textbf{56} (2001), no. 1, 185--186.

\bibitem{Spir2} V.P. Spiridonov, \textit{Theta hypergeometric integrals}, Algebra i Analiz \textbf{15} (2003), 161--215 (St Petersburg Math J. \textbf{15} (2004), 929--967).


\end{thebibliography}
\end{document}